\documentclass[12pt]{article}
\usepackage{amssymb,amsthm,amsmath,amsfonts,latexsym,tikz,hyperref,youngtab,mathdots,url}
\usepackage[hmargin=.8in,vmargin=.5in]{geometry}

\newcommand{\ben}{\begin{enumerate}}
\newcommand{\een}{\end{enumerate}}
\newcommand{\ble}{\begin{lem}}
\newcommand{\ele}{\end{lem}}
\newcommand{\bth}{\begin{thm}}
\renewcommand{\eth}{\end{thm}}
\newcommand{\bpr}{\begin{prop}}
\newcommand{\epr}{\end{prop}}
\newcommand{\bco}{\begin{cor}}
\newcommand{\eco}{\end{cor}}
\newcommand{\bcon}{\begin{conj}}
\newcommand{\econ}{\end{conj}}
\newcommand{\bde}{\begin{defn}}
\newcommand{\ede}{\end{defn}}
\newcommand{\bex}{\begin{exa}}
\newcommand{\eex}{\end{exa}}
\newcommand{\barr}{\begin{array}}
\newcommand{\earr}{\end{array}}
\newcommand{\btab}{\begin{tabular}}
\newcommand{\etab}{\end{tabular}}
\newcommand{\beq}{\begin{equation}}
\newcommand{\eeq}{\end{equation}}
\newcommand{\bea}{\begin{eqnarray*}}
\newcommand{\eea}{\end{eqnarray*}}
\newcommand{\bal}{\begin{align*}}
\newcommand{\bce}{\begin{center}}
\newcommand{\ece}{\end{center}}
\newcommand{\bpi}{\begin{picture}}
\newcommand{\epi}{\end{picture}}
\newcommand{\bpp}{\begin{picture}}
\newcommand{\epp}{\end{picture}}
\newcommand{\bfi}{\begin{figure} \begin{center}}
\newcommand{\efi}{\end{center} \end{figure}}
\newcommand{\bprf}{\begin{proof}}
\newcommand{\eprf}{\end{proof}\medskip}
\newcommand{\capt}{\caption}
\newcommand{\bsl}{\begin{slide}{}}
\newcommand{\esl}{\end{slide}}
\newcommand{\bfr}{\begin{frame}}
\newcommand{\efr}{\end{frame}}

\newtheorem{thm}{Theorem}[section]
\newtheorem{prop}[thm]{Proposition}
\newtheorem{cor}[thm]{Corollary}
\newtheorem{lem}[thm]{Lemma}
\newtheorem{conj}[thm]{Conjecture}
\newtheorem{exa}[thm]{Example}

\newcommand{\lf}{\lfloor}
\newcommand{\rf}{\rfloor}

\newcommand{\binst}[2]{\left\{ #1 \atop #2 \right\}}
\newcommand{\gaus}[2]{\left[  \left.#1 \atop #2 \right.  \right]}

\newcommand{\bbP}{{\mathbb P}}
\newcommand{\ree}[1]{(\ref{#1})}
\newcommand{\qmq}[1]{\quad\mbox{#1}\quad}
\newcommand{\De}{\Delta}
\newcommand{\ze}{\zeta}
\newcommand{\eqqed}[1]{$\rule{1ex}{0ex}\hfill{\dil#1}\hfill\qed$}
\newcommand{\dil}{\displaystyle}
\newcommand{\vs}[1]{\vspace{#1}}
\newcommand{\cC}{{\cal C}}
\newcommand{\cL}{{\cal L}}
\DeclareMathOperator{\wt}{wt}
\newcommand{\ep}{\epsilon}
\newcommand{\hqed}{\hfill \qed}
\newcommand{\case}[4]{\left\{\barr{ll}#1&\mbox{#2}\\#3&\mbox{#4}\earr\right.}
\newcommand{\de}{\delta}
\newcommand{\bbZ}{{\mathbb Z}}
\newcommand{\hqedm}{\hfill \qed \medskip}
\newcommand{\bb}{{\bf b}}
\newcommand{\bF}{{\bf F}}
\newcommand{\ka}{\kappa}
\DeclareMathOperator{\Mod}{mod}

\begin{document}
\pagestyle{plain}

\title{Generalized Fibonacci polynomials and Fibonomial coefficients
}
\author{
Tewodros Amdeberhan\\[-5pt]
\small Department of Mathematics, Tulane University,\\[-5pt]
\small New Orleans, LA 70118, USA, {\tt tamdeber@tulane.edu}\\
\\
Xi Chen\\[-5pt]
\small School of Mathematical Sciences, Dalian University of Technology,\\[-5pt]
\small Dalian City, Liaoning Province, 116024, P. R. China, {\tt xichen.dut@gmail.com} 
\thanks{Research partially supported by a grant from the China Scholarship Council}\\
\\
Victor H. Moll\\[-5pt]
\small Department of Mathematics, Tulane University,\\[-5pt]
\small New Orleans, LA 70118, USA, {\tt vhm@tulane.edu}
\thanks{Research partially supported by the National Science Foundation NSF-DMS 1112656}\\
\\
Bruce E. Sagan\\[-5pt]
\small Department of Mathematics, Michigan State University,\\[-5pt]
\small East Lansing, MI 48824, USA, {\tt sagan@math.msu.edu}
}

\date{\today\\[10pt]
	\begin{flushleft}
	\small Key Words:  binomial theorem, Catalan number, Dodgson condensation, Euler-Cassini identity, Fibonacci number, Fibonomial coefficient, Lucas number $q$-analogue, valuation
	                                       \\[5pt]
	\small AMS subject classification (2000): 
	Primary 05A10;
	Secondary 11B39, 11B65.\\[5pt]
	\small Running title: Fibonacci polynomials and Fibonomial coefficients
	\end{flushleft}}

\maketitle

\begin{abstract}
The focus of this paper is the study of generalized Fibonacci polynomials and Fibonomial coefficients.  The former are polynomials $\{n\}$ in variables $s,t$ given by $\{0\}=0$, $\{1\}=1$, and $\{n\}=s\{n-1\}+t\{n-2\}$ for $n\ge2$.  The latter are defined by 
$\binst{n}{k}=\{n\}!/(\{k\}!\{n-k\}!)$ where $\{n\}!=\{1\}\{2\}\dots\{n\}$.  These quotients are also polynomials in $s,t$ and specializations give the ordinary binomial coefficients, the Fibonomial coefficients, and the $q$-binomial coefficients.  We present some of their fundamental properties, including a more general recursion for $\{n\}$, an analogue of the binomial theorem, a new proof of the Euler-Cassini identity in this setting with applications to estimation of tails of series, and valuations when $s$ and $t$ take on integral values.  We also study a corresponding analogue of the Catalan numbers.  Conjectures and open problems are scattered throughout the paper.
\end{abstract}

%
%

\section{Introduction}
\label{sec-intro}

We will be studying generalized Fibonacci polynomials and generalized Fibonomial coefficients.   Throughout this work, $\bbP$ will stand for the positive integers.  The \emph{Fibonacci numbers} $F_{n}$ are defined by $F_0=0$, $F_1=1$ and, for $n\ge2$,
$$
F_{n} = F_{n-1} + F_{n-2}.
$$
The \emph{Lucas numbers}
$L_{n}$ are defined by the same recurrence, with the initial conditions
$L_{0}=2$ and $L_{1}=1$. The reader will find an 
introduction to these well-studied sequences in the books by Koshy~\cite{kos:fln} and Moll~\cite{mol:nfs}.

One  generalization of these numbers which has received much attention is the sequence of \emph{Fibonacci polynomials}
\beq
\label{Fn(x)}
F_{n}(x) = x F_{n-1}(x) + F_{n-2}(x), \quad n \geq 2,
\eeq
with initial conditions $F_{0}(x) = 0, \, F_{1}(x) = 1$. 
The \emph{generalized Fibonacci polynomials} which  we will consider depend on 
two variables 
$s, t$ and are  defined by  $\{ 0 \}_{s,t} = 0$, $\{ 1 \}_{s,t} = 1$ and, for $n\ge2$,
\beq
\label{nst-rr}
\{ n \}_{s,t} = 
s \{n-1 \}_{s,t} + t \{ n-2 \}_{s,t}.
\eeq
Here and with other quantities depending on $s$ and $t$, we will often drop the subscripts if they are clear from context.  For example, we have
$$
\{ 2 \}= s, \quad \{ 3 \} = s^{2}+t, \quad \{ 4 \}= s^{3} + 2st, 
\quad \{ 5 \} = s^{4} + 3s^{2}t+t^{2}.
$$
When $s$ and $t$ are integers, these sequences were first studied by Lucas in a series of papers~\cite{luc:tfn1,luc:tfn2,luc:tfn3} and then forgotten.  Nearly 100 years later,  Hoggatt and  Long \cite{hl:dpg} rediscovered them, this time considering $s$ and $t$ as variables.  But they have received considerably less attention than the one variable family in~\ree{Fn(x)}, although some of their properties are the same because of the relation
$$
\{n\}_{s,t}=t^{(n-1)/2} F_n\left(\frac{s}{\sqrt{t}}\right).
$$
Part of the purpose of the present work  is to rectify this neglect.

Our notation is chosen to reflect two important specializations of this sequence (other than the one $s=t=1$ already mentioned). In particular, if $s=2$ and $t=-1$ then $\{n\}=n$.  And if $s=q+1$ and $t=-q$ then 
\beq
\label{[n]}
\{n\}=1+q+q^2+\dots+q^{n-1}=[n]_q,
\eeq
the standard $q$-analogue of $n$.

There is a  corresponding extension of the Lucas numbers, the \emph{generalized Lucas polynomials}, defined by 
$$
\langle n \rangle_{s,t} = s \langle n-1 \rangle_{s,t}  + 
t \langle  n-2 \rangle_{s,t} , \quad n\geq 2
$$
together with the initial conditions $\langle  0 \rangle_{s,t} = 2$ 
and $\langle 1 \rangle_{s,t}  = s$.  Here is a list of the first few polynomials
\begin{equation*}
\langle 2 \rangle_{s,t} = s^{2}+2t, \quad \langle 3 \rangle_{s,t} = s^{3}+3st, \quad
\langle  4 \rangle_{s,t} = s^{4} + 4s^{2}t+2t^{2}, 
\quad \langle 5 \rangle_{s,t} = s^{5} + 5s^{3}t+5st^{2}.
\end{equation*}
Of course, when $s=t=1$ these reduce to the ordinary Lucas numbers.

One can find algebraic expressions for these polynomials using standard techniques from the theory of recursively defined sequences.  In particular, the characteristic polynomial of the recurrence is $z^2-sz-t$ whose roots are
\beq
\label{XY}
X = \frac{s +\sqrt{s^{2}+4t}}{2}
\qmq{and}
Y = \frac{s -\sqrt{s^{2}+4t}}{2}.
\eeq
We will often abbreviate $ \sqrt{s^{2}+4t}=\De_{s,t}$. 
The following relations between $s,t$ and $X,Y$ will be useful
\beq
\label{stXY}
s=X+Y,\quad t = -XY,\quad \De= X-Y.
\eeq

Any solution of~\ree{nst-rr} can be expressed as $c X^n+ d Y^n$ for constants $c,d$ depending on the initial conditions.  Computing the constants in the two cases of interest to us gives the following analogue of Binet's formula for the Fibonacci numbers.
\bpr
\label{nXY}
For $n\ge0$ we have

\medskip

\eqqed{
\{ n \} = \frac{X^{n} - Y^{n}}{X-Y}  \qmq{and} \langle n \rangle = X^{n} + Y^{n}.
}
\epr

Combining this result with equation~\ree{[n]}, we see that there is another relation between $\{n\}$ and $[n]$, namely
\beq
\label{[n]_{X/Y}}
\{n\}_{s,t}=Y^{n-1} [n]_{X/Y}.
\eeq
This will be useful in the sequel.

\bfi
\begin{tikzpicture}
\draw (0,0) grid (3,1);
\fill (.5,.5) circle (.1);
\fill (1.5,.5) circle (.1);
\fill (2.5,.5) circle (.1);
\end{tikzpicture}
\qquad
\begin{tikzpicture}
\draw (0,0) grid (3,1);
\draw (.5,.5)--(1.5,.5);
\fill (.5,.5) circle (.1);
\fill (1.5,.5) circle (.1);
\fill (2.5,.5) circle (.1);
\end{tikzpicture}
\qquad
\begin{tikzpicture}
\draw (0,0) grid (3,1);
\draw (1.5,.5)--(2.5,.5);
\fill (.5,.5) circle (.1);
\fill (1.5,.5) circle (.1);
\fill (2.5,.5) circle (.1);
\end{tikzpicture}

\vs{10pt}

\begin{tikzpicture}
\draw (0,0) circle (.4);
\draw (0,0) circle (1);
\draw (30:.4) -- (30:1)  (150:.4) -- (150:1) (270:.4) -- (270:1);
\fill (90:.7) circle (.1);
\fill (210:.7) circle (.1);
\fill (330:.7) circle (.1);
\end{tikzpicture}
\qquad
\begin{tikzpicture}
\draw (0,0) circle (.4);
\draw (0,0) circle (1);
\draw (30:.4) -- (30:1)  (150:.4) -- (150:1) (270:.4) -- (270:1);
\fill (90:.7) circle (.1);
\fill (210:.7) circle (.1);
\fill (330:.7) circle (.1);
\draw (210:.7) arc (210:330:.7);
\end{tikzpicture}
\qquad
\begin{tikzpicture}
\draw (0,0) circle (.4);
\draw (0,0) circle (1);
\draw (30:.4) -- (30:1)  (150:.4) -- (150:1) (270:.4) -- (270:1);
\fill (90:.7) circle (.1);
\fill (210:.7) circle (.1);
\fill (330:.7) circle (.1);
\draw (90:.7) arc (90:-30:.7);
\end{tikzpicture}
\qquad
\begin{tikzpicture}
\draw (0,0) circle (.4);
\draw (0,0) circle (1);
\draw (30:.4) -- (30:1)  (150:.4) -- (150:1) (270:.4) -- (270:1);
\fill (90:.7) circle (.1);
\fill (210:.7) circle (.1);
\fill (330:.7) circle (.1);
\draw (90:.7) arc (90:210:.7);
\end{tikzpicture}
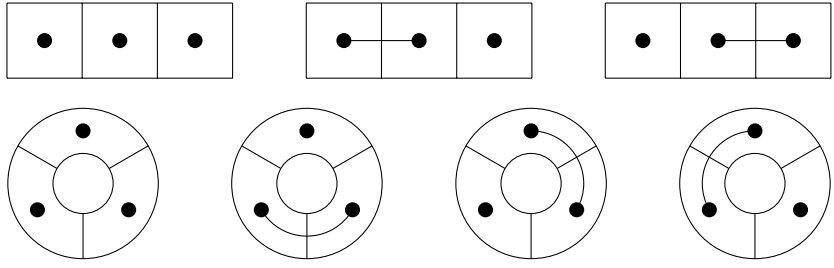
\capt{ Linear and circular tilings}
\label{fig-tile}
\efi

In addition to this algebraic approach to our polynomials, there is a  combinatorial interpretation derived from  the standard interpretation of $F_n$ via tiling.  A \emph{linear tiling}, $T$,  of a row of squares is a covering of the squares with dominos (which cover two squares) and monominos (which cover one square).  We let
$$
\cL_n=\{ T\ :\ \text{$T$ a linear tiling of a row of $n$ squares}\}.
$$
The  three tilings in the first row of  Figure~\ref{fig-tile} are the elements of $\cL_3$.  We will also consider \emph{circular tilings} where the (deformed) squares are arranged in a circle.  We will use the notation $\cC_n$ for the set of circular tilings of $n$ squares.  So the set of  tilings in the bottom row of  Figure~\ref{fig-tile} is $\cC_3$.  For any type of nonempty tiling, $T$, we define its \emph{weight} to be
$$
\wt T = s^{\text{\# of monominos in $T$}} t^{\text{\# of dominos in $T$}}.
$$
We give the empty tiling $\ep$ of zero boxes the weight $\wt\ep=1$ if it is being considered as a linear tiling or $\wt\ep=2$ if it is being considered as a circular tiling.  The following proposition is immediate from the definitions of weight and of our generalized polynomials.
\bpr
\label{wt}
For $n\ge 0$ we  have

\medskip

\eqqed{
\{n+1\} =\sum_{T\in \cL_n} \wt T \qmq{and} \langle n \rangle = \sum_{T\in \cC_n} \wt T.
}
\epr

We are now in a position to define corresponding generalized binomial coefficients.
Given a sequence $\mathbf{a} = \{ a_{n} \}$ of nonzero real numbers, it is natural to define the
$\mathbf{a}$-factorials  by 
$$
n!_{\mathbf{a}}  = \prod_{i=1}^{n} a_{i},
$$
\noindent
and the $\mathbf{a}$-binomial coefficients by 
$$
\binom{n}{k}_{\mathbf{a}}  = \frac{n!_{\mathbf{a}}}
{k!_{\mathbf{a}} \, (n-k)!_{\mathbf{a}}}. 
$$
The classical example comes from $a_{n} = n$  with the standard $q$-analogue being obtained when $a_n=[n]_q$.  Here, we will consider
$a_{n} = \{ n \}_{s,t}$ and $a_n=\langle n \rangle_{s,t}$. To simplify the notation, these 
are written, respectively, as 
$$
\{ n \}_{s,t}! \quad \text{ and } \quad \langle n \rangle_{s,t}! 
$$
for factorials, and 
$$
\left\{ n \atop k \right\}_{s,t} \quad \text{ and } \quad
\left\langle n \atop k \right\rangle_{s,t}
$$
for binomial coefficients.  The product $\{n\}!$ is called a \emph{generalized fibotorial} and $\langle n\rangle !$ is a \emph{generalized lucatorial}.  
The binomial coefficients $\left\{ n \atop k  \right\}$ and 
$\left\langle n \atop k \right\rangle$ 
are called generalized \emph{Fibonomial} and \emph{Lucanomial} coefficients, respectively.

We can relate the generalized Fibonomials to the $q$-binomial coefficients algebraically.  Indeed, it follows easily from~\eqref{[n]_{X/Y}} that 
\beq
\label{gau}
\binst{n}{k}_{s,t}= Y^{k(n-k)} \gaus{n}{k}_{X/Y}.
\eeq
There is also a simple combinatorial interpretation of $\binst{n}{k}$ which was given by Savage and Sagan~\cite{ss:cib} using tilings of a $k\times(n-k)$ rectangle containing a partition.  But we will not use this description, and instead we refer the reader to their paper for the details.  It would be very interesting to give combinatorial proofs for some of the results we give for Fibonomials.

The rest of this paper is structured as follows.  In the next section we present some of the fundamental properties of $\{n\}$ and $\langle n \rangle$ which will be useful in the rest of the paper.  Section~\ref{sec-arithmetic} is devoted to valuations.  In particular, we provide a complete description of the $2$-adic valuation of $\{n\}$ and $\{n\}!$ for arbitrary integers $s,t$.  The Euler-Cassini identity is the focus of Section~\ref{sec-EC}.  We give a new proof of the diagonal case of this equality using Dodgson condensation.  We also use this identity to give estimates for the tails of Fibonacci analogues of the series for the Riemann zeta function.  In Section~\ref{sec-fpb} we return to the generalized Fibonomial coefficients and consider their recursions and analogues of the binomial theorem.  We end by studying an $s,t$-version of the Catalan numbers $C_n$  and proving analogues of the theorem giving the $2$-adic valuation of $C_n$.  
Various open problems and conjectures are mentioned.

%
%

\section{Fundamental properties of $\{n\}$ and $\langle n\rangle$}
\label{sec-fpn}

In this section we collect some of the important properties of 
$\{ n \}_{s,t}$ and $\langle n \rangle_{s,t}$ to be used in the sequel.  Most can be proved using either the algebraic descriptions in terms of $X$ and $Y$, or the combinatorial interpretations, or both.  Often the demonstrations do not differ significantly from ones already in the literature, so we will sometimes omit the proofs and just supply references where they can be found.

We begin with the expansion of our polynomials into monomials.
\bpr
\label{{n}-st}
The polynomials $\{ n \}$ and $\langle n \rangle$ are given by 
\beq
\{ n \} = \sum_{k \geq 0} \binom{n-k-1}{k} s^{n-2k-1}t^{k}
\label{sum-1a}
\eeq
and 
$$
\langle n \rangle = \sum_{k \geq 0} \frac{n}{n-k} \binom{n-k}{k}
s^{n-2k} t^{k}.
$$
\epr
\begin{proof}
We just prove the first identity as the second is similar.  By Proposition~\ref{wt}, it suffices to show that the $k$th term of the sum is the sum of the weights of all linear tilings of $n-1$ squares with $k$ dominos.  If there are $k$ dominos, then there must be $n-2k-1$ monominos and so this accounts for the monomial $s^{n-2k-1}t^{k}$.  To count the number of arrangements of these tiles, number the squares $1,\dots,n-1$ from left to right.  Then picking the places for the left endpoints of the dominos is equivalent to picking $k$ numbers from $1,\dots,n-2$ with no two consecutive.  The number of ways of doing this is $\binom{n-k-1}{k}$, and once the dominos are placed, there is no further  choice for distributing the monominos.  This finishes the proof.
\end{proof}

The next result is a useful generalization of the defining recurrence for the polynomials $\{ n \}$. 

\bth
\label{thm-recu1}
For $m\ge1$ and $n\ge0$ we have
$$
\{ m + n \} = \{ m \} \{ n+1 \} + t \{ m-1 \} \{ n \}.
$$
\eth
\begin{proof}
This can be given a combinatorial proof (see~\cite{ss:cib}), but we will indicate an algebraic one to illustrate the method.  First one uses~\eqref{stXY} and Proposition~\ref{nXY} to convert the equality into an equivalent statement about polynomials in $X$ and $Y$.  This statement can then be easily verified by algebraic manipulations.
\end{proof}

We now note two identities relating $\{n\}$ and $\langle n\rangle$.
\bpr[\cite{ss:cib}]
For $n \geq 1$ we have
$$
\langle n \rangle = \{ n +1  \} + t \{ n - 1 \}.
$$
And for $m,n\ge0$ we have

\medskip

\eqqed{
\{ m+ n \} = \frac{ \langle m \rangle \{ n \} + 
\{ m \} \langle n \rangle }{2}.
}
\epr

The next result can be used to show that many divisibility properties carry over directly from the integers $n$ to the polynomials $\{n\}$. 

\bpr[\cite{hl:dpg}]
\label{gcd}
For $m,n\ge1$ we have
$$
\gcd( \{ m \}, \{n \} ) = \{ \gcd(m,n) \}.
$$
\noindent
Equivalently, $m$ divides $n$ if and only if $\{m\}$ divides $\{n\}$.\hqed
\epr

Using standard techniques, one can convert the defining recursion for $\{n\}$ into a generating function.

\bpr
\label{gf}
The generating function of the polynomials $\{ n \}$ is given by 

\medskip

\eqqed{
 \sum_{n = 0}^{\infty} \{ n \} z^{n} = \frac{z}{1 - sz - tz^{2}}.
}
\epr

As an application of this result, we will derive a generalization of the following well-known identity for Fibonacci numbers
$$
\sum_{n=0}^{\infty} \frac{F_{n}}{2^{n+1}} = 1
$$
which is the special case $s=t=1$ and $z=1/2$ of the above proposition.

\bco
\label{nice-1}
For $s, t\in \bbP$ we have
$$
\sum_{n =0}^{\infty} \frac{t \{ n \}_{s,t}}{(s+t)^{n+1}} = \frac{1}{s+t-1}.
$$
\eco
\begin{proof}
We will give both algebraic and combinatorial proofs.
The former is obtained by setting $z = 1/(s+t)$ in the generating 
function of Proposition~\ref{gf}.  We must make sure that this substitution is analytically valid in that  $1/(s+t)$ is smaller than the radius of convergence of the power series which is $\min\{1/|X|,1/|Y|\}$.  But this is a routine check using equation~\eqref{XY}.  Once the substitution is made, simple algebraic manipulations complete the demonstration.

For the combinatorial proof,
consider an infinite row of squares numbered left to right by the positive  integers.  Suppose each square can be colored with one of
$s$ shades of white and $t$ shades of black. 
Let $Z$ be the random variable which returns the box number at the end of the first odd-length block of boxes all of the same black shade.   For $n \in \bbP$,  the event $Z = n$  is equivalent to having box $n$ painted with one of the 
shades of black, box $n+1$ painted with any of the remaining colors, and all blocks of a black shade among the first $n-1$ squares being of even length.    So there are $t$ choices for the color of box $n$ and $s+t-1$ choices for the color of box $n+1$.   Each coloring of the first $n-1$ squares gives rise to a tiling where each white box is replaced by a monomino and a block of $2k$ boxes of the same black shade is replaced by $k$ dominoes.  Also, the weight of the tiling is just the number of colorings mapping to it.  Thus, by Proposition~\ref{wt}, the number of colorings for the first $n-1$ boxes is $\{n\}$.  Hence 
$$
P(Z=n) = \frac{t (s+t-1) \, \{ n \}_{s,t}}{(s+t)^{n+1}}.
$$
Summing these probabilities finishes the proof.
\end{proof}

We end this section by exploring the binomial transformation of the sequence $\{n\}$, $n\ge0$.  Interestingly, doing so involves a change of variables from $s,t$ to $s+2,t-s-1$.  We then use the transform to prove a well-known identity for Fibonacci numbers.
\bpr
For $n\ge0$ we have
$$
\sum_{k=0}^{n} \binom{n}{k} \{ k \}_{s,t} = \{ n \}_{s+2,t-s-1}.
$$
In particular, for $s=t=1$,
$$
\sum_{k=0}^{n} \binom{n}{k}F_{k} =  F_{2n}.
$$
\epr
\begin{proof}
Using Proposition~\ref{nXY}, we have the exponential generating function
$$
\sum_{k=0}^{\infty}\{k\}_{s,t}\frac{z^k}{k!}=\frac{e^{Xz}-e^{Yz}}{X-Y}.
$$
Note that
$$
X+1, Y+1=\frac{(s+2)\pm\sqrt{s^2+4t}}2=\frac{(s+2)\pm\sqrt{(s+2)^2+4(t-s-1)}}2.
$$
Putting everything together and using the product rule for exponential generating functions gives
$$
\sum_{n=0}^{\infty}\frac{z^n}{n!}\sum_{k=0}^n\binom{n}k\{k\}_{s,t}=
e^z\frac{e^{Xz}-e^{Yz}}{X-Y}=\frac{e^{(X+1)z}-e^{(Y+1)z}}{X-Y}=
\sum_{n\geq0}\{n\}_{s+2,t-s-1}\frac{z^n}{n!}.
$$ 
Extracting the coefficients of $z^n/n!$ completes the proof of the first equation.

For the second, from what we have just proved it suffices to show that $\{n \}_{3,-1} = \{ 2 n \}_{1,1} $.  But
$$
 \{ 2 n \}_{1,1}
=\frac{1}{\sqrt{5}}\left[ \left( \frac{1+\sqrt{5}}{2} \right)^{2n} -   \left( \frac{1-\sqrt{5}}{2} \right)^{2n} \right]
=\frac{1}{\sqrt{5}}\left[ \left( \frac{3+\sqrt{5}}{2} \right)^{n} -   \left( \frac{3-\sqrt{5}}{2} \right)^{n} \right]
=\{n\}_{3,-1}
$$
and we are done.
\end{proof}

%
%

\section{Arithmetic properties}
\label{sec-arithmetic}

We will be  concerned with the \emph{$d$-adic valuation} function
$$
\nu_d(n) = 
\begin{cases}
\text{the highest power of $d$ dividing $n$}&\text{if $n\neq 0$,}\\
\infty & \text{if $n=0$.}
\end{cases}
$$
If the subscript is missing, then it is assumed that $d=2$.
A fact about valuations which we will use repeatedly is that if $\nu_d(m)\neq\nu_d(n)$ then 
\beq
\label{min}
\nu_d(m+n)=\min\{\nu_d(m),\nu_d(n)\}.
\eeq
Our primary goal will be to characterize $\nu_2(\{n\}_{s.t})$ for all possible integers $s,t$.  This will then be used in Section~\ref{sec-Catalan} to give analogues of a well-known theorem about the $2$-adic valuation of the Catalan numbers.  We will end the section with an indication of what can be said for other moduli.

We will now characterize $\nu(\{n\})=\nu_2(\{n\}_{s,t})$ for all integral $s,t$ as well as $\nu(\{n\}!)$.  We first consider the case when both $s,t$ are odd.   If $S$ is a set of integers then we will  have much use for the indicator function
$$
\de_S(k)=
\begin{cases}
1&\text{if $k\in S$,}\\
0&\text{if $k\not\in S$.}
\end{cases}
$$
In this context, we will let $E$ and $O$ stand for the even and odd integers, respectively.
\ble
\label{lem-sodd-todd}
Let $s$ and $t$ be odd.  We have $\nu(\{n\})=0$ whenever $n=3k+1$ or $3k+2$.  If $n=3k$ then
$$
\nu(\{3k\})=\case{1+\de_E(k)(\nu(k\{6\})-2)}{if $t\equiv 1 \pmod{4}$,}{\nu(k\{3\})}{if $t\equiv 3 \pmod{4}$.\rule{0pt}{20pt}}
$$
\ele
\bprf
Our proof will be by induction on $n$ where the base cases are easy to verify.  
From the recursion 
\beq
\label{odd1}
\{n\}=s\{n-1\}+t\{n-2\}
\eeq
and fact that $s$ and $t$ are odd, it is clear that both $\{3k+1\}$ and $\{3k+2\}$ are odd while $\{3k\}$ is even.  This finishes the cases when $n=3k+1$ or $n=3k+2$. The demonstrations when $n$ is divisible by $3$ are similar for both possible residues of $t$, so we will only present $t\equiv 3 \pmod{4}$.

Suppose now that $n=6k+3$ for some integer $k$.   Using the recursion in Theorem~\ref{thm-recu1} we have
$$
\{6k+3\}=\{3\}\{6k+1\}+t\{2\}\{6k\}.
$$
By hypothesis and induction we know that $\{6k+1\}$, $\{2\}$, and $t$ are odd.  Furthermore, by induction again,
$\nu(\{6k\})=1+\nu(k\{3\})>\nu(\{3\})$.  So using~\ree{min} on the previous displayed equation gives $\nu(\{6k+3\})=\nu(\{3\})=\nu((n/3)\{3\})$ since $n/3$ is odd.  This is the desired conclusion.

For the final case, let $n=6k+6$ for some $k$.  Using Theorem~\ref{thm-recu1} repeatedly we obtain
$$
\{6k+6\}=\{3k+4\}\{3k+3\}+t\{3k+2\}\{3k+3\}
=\{3k+3\}((s^3+3st)\{3k+1\}+(s^2t+2t^2)\{3k\}).
$$
As before, we can ignore $\{3k+1\}$ and factors of $s$ or $t$ since they are odd.  Since $s$ is odd, $s^2\equiv 1 \pmod{4}$.  It follows that  $\nu(s^2+3t)=1$ while $\nu(\{3\})=\nu(s^2+t)\ge2$.  Applying~\ree{min} and induction to the previous displayed equation gives 
$$
\nu(\{6k+6\})=\nu(\{3k+3\})+1=\nu((k+1)\{3\})+1=\nu((2k+2)\{3\})=\nu((n/3)\{3\})
$$
which is, again, what we want.
\eprf

We can use this lemma to calculate the $2$-adic valuation of the corresponding  factorials.
\bco
\label{cor-sodd-todd}
Let $s$ and $t$ be odd.   We have
$$
\nu(\{n\}!)=
\begin{cases}
\nu(\lf n/3 \rf !) + \lf n/6 \rf \nu(\{6\})+\de_O(\lf n/3)\rf) & \text{if $t\equiv1\pmod{4}$,}\\[10pt]
\nu(\lf n/3 \rf !)+\lf n/3 \rf \nu(\{3\}) & \text{if $t\equiv3\pmod{4}$,}
\end{cases}
$$
where $\lf\cdot\rf$ is the floor function.\hfill\qed
\eco
\bprf
Again, we will only provide a proof when $t\equiv 3 \pmod{4}$ as the other congruence class is similar.  Write $n=3k+r$ where $0\le r<3$.  Using Lemma~\ref{lem-sodd-todd} we have
$$
\nu(\{n\}!)=\sum_{i=1}^n \nu(\{i\})=\sum_{i=1}^k \nu(\{3i\})=\sum_{i=1}^k \nu(i\{3\})=\nu(k!)+k\nu(\{3\})
$$
and using the fact that $k=\lf n/3\rf$ finishes the demonstration.
\eprf

Now we turn to the cases  when $s$ and $t$ are of opposite parity.  If $s$ is odd and $t$ is even then a simple induction shows that $\{n\}$ is always odd for $n\ge1$.  The reverse case is more interesting.

\ble
\label{lem-seven-todd}
Let $s$ be even and $t$ be odd.  We have
$$
\nu(\{n\})=\case{\nu(sn/2)}{if $n$ is even,}{0}{if $n$ is odd.}
$$
\ele
\bprf
This proof is much like the one for Lemma~\ref{lem-sodd-todd} and so we will content ourselves with stating the main equation for the induction step on even integers
\beq
\label{eq-2n}
\{2n\}=\{n\}(s\{n\}+2t\{n-1\}).
\eeq
The reader can easily fill in the details.
\eprf

The proof of the following corollary is much like that of Corollary~\ref{cor-sodd-todd}  and so is omitted.
\bco
\label{cor-seven-todd}
Let $s$ be even and $t$ be odd.  We have

\smallskip

\eqqed{
\nu(\{n\}!)=\nu(n!)+\lf n/2\rf\nu(s/2).
}

\eco

Finally, we have the case where both parameters are even.  To describe the $2$-adic valuations we will rely on a recursion.

\begin{lem}
Let $s,t$ be odd.  We have
\beq
\label{eq-base}
\nu(\{n\}_{2^a s, 2t})=
\begin{cases}
\lf n/2 \rf + \de_{4\bbZ}(n)\left[\nu(n \{4\}_{2s,2t}/4)-2\right]	&\text{if $a=1$,}\\
\lf n/2 \rf + \de_E(n)[\nu(n)+a-2]						&\text{if $a\ge2$.}
\end{cases}
\eeq
Now suppose $s,t$ are arbitrary integers.  We have
\beq
\label{eq-2s-4t}
\nu(\{n\}_{2s,4t})=n-1+\nu(\{n\}_{s,t}).
\eeq
This completely determines the $2$-adic valuations of $\{n\}$ where both subscripts are even.
\end{lem}
\bprf To prove~\ree{eq-base} the usual ideas come into play.  The equations which are used for the induction are the defining recursion, equation~\ree{eq-2n}, and
$$
\{8k+4\}_{s,t}=\{4\}_{s,t}\{8k+1\}_{s,t}+t\{3\}_{s,t}\{8k\}_{s,t}.
$$

The proof of  equation~\ree{eq-2s-4t} is very simple.  In fact, a straightforward induction on $n$ shows that $\{n\}_{2s,4t}=2^{n-1} \{n\}_{s,t}$ which implies the desired result.

For the last statement, by repeated use of equation~\ree{eq-2s-4t}, one can reduce finding $\nu(\{n\})$ to finding $\nu(\{n\}_{s.t})$ where either at least one of $s, t$ is odd or both are even and $t$ is twice an odd number.  In the former case, the computation is finished by one of our former results.  In the latter case, one can use equation~\ree{eq-base} to complete the evaluation.
\eprf

Because of the recursive nature of these valuations, the corresponding formulas for $\nu(\{n\}!)$ are complicated and too messy to be of real interest.  On the other hand, we do not wish to give the impression that one can only say interesting things for the modulus $d=2$.  So our last result in this section will be for arbitrary $d$.

\bpr
Consider any positive integer $d\ge2$.  We have, for any $n\ge1$,
$$
\nu_d(\{n\}_{d,-1}) = \de_E(n) \nu_d(dn/2).
$$
\epr
\bprf
First, consider the case where $d$ is a  prime.  Using the defining recursion for $\{n\}=\{n\}_{d,-1}$ one easily sees that $\{n\}$ is divisible by $d$ if and only if $n$ is even.   This completes the $n$ odd case.  For even integers, letting $s=d$ and $t=-1$ in 
equation~\ree{sum-1a} and re-indexing gives
$$
(-1)^{n-1} \{2n\}=\sum_{k\ge0} \binom{n+k}{2k+1} d^{2k+1}(-1)^k
=dn + dn\sum_{k\ge1}\frac{c_k d^{2k}}{(2k+1)!}
$$
where the $c_k$ are integers because $n$ divides $(n+k)(n+k-1)\dots(n-k)$. 
So, by equation~\ree{min}, it suffices to show that the $d$-adic valuation of every term in the sum over positive $k$ is at least $1$.  Since $d$ is prime, we can use Legendre's well-known formula
$\nu_d(n!) = \sum_{i\ge1} \lf n/d^i \rf$ to show that
$$
\nu_d((2k+1)!)\le
\begin{cases}
\dil\frac{2k}{d-1} &\text{if $d\ge3$,}\\[10pt]
2k-1  &\text{if $d=2$.}
\end{cases}
$$
From this, it is easy to verify that $\nu_d(d^{2k}/(2k+1)!)\ge1$ which completes the case when $d$ is prime.

To finish the proof, note first that the only place where we used the fact that $d$ was prime was in deriving the upper bounds  on $\nu_d((2k+1)!)$.  But these will still hold when $d$ is a prime power, and may even become sharper.  Finally, for general $d$ we just use the fact that  if $p$ and $q$ are relatively prime then 
$\nu_{pq}(n)=\min\{ \nu_p(n),\nu_q(n)\}$ for any integer $n$.
\eprf

We conjecture that  the roles of the modulus $d$ and the parameter $s$ in the previous proposition can be decoupled.
\bcon
Suppose $s\ge2$ is an integer and $d\ge3$ is an  odd integer.  There exist positive integers $s^*, d^*$ depending only on $s,d$ such that
$d^*\le d$ and
$$
\nu_d(\{n\}_{s,-1}) = \de_{d^*\bbZ}(n) \nu_d(s^*n/d^*).
$$
\econ 
\noindent Of course, it would be desirable to have a way of computing $s^*$ and $d^*$ from $s$ and $d$ rather than just an existential proof.  The following table lists pairs $(s^*,d^*)$ for  small values of $s$ and $d$ indexing the rows and columns, respectively.   Note that any two values of $s^*$ with the same valuation modulo $d$ will yield the same result on the right-hand side of the equation in the conjecture.  Also note that  if the first positive integer $n$ with $\nu_d(\{n\}_{s,-1})\neq 0$ is $n=d$ then we have a choice as to whether to let $d^*=1$ or $d^*=d$.   
$$
\barr{l|cccc}
s\setminus d& 
       	3&	5&	7&	9\\
\hline
2    &  (1,1)&	(1,1)&	(1,1)&	(1,1)\\
3    & (3,2)&	(1,1)&	(7,4)&	(9,6)\\
4    & (1,1)&	(5,3)&	(7,4)&  (1,1)\\
5    & (1,1)&	(5,2)&	(1,1)&	(1,1)
\earr
$$

%
%

\section{The Euler-Cassini identity}
\label{sec-EC}

In this section we will consider various results related to the famous Euler-Cassini identity for Fibonacci numbers.  We will first recall a version of this equation proved by Cigler~\cite{cig:qfp} for a $q$-analogue of our generalized Fibonacci polynomials.  We then show how Dodgson condensation~\cite{dod:cd} can be used to prove a particular case of this identity.  Finally, we use a slightly more general form of Euler-Cassini to give estimates for the tails of certain infinite  series with terms involving the polynomials $\{n\}$ evaluated at various integers.

The $q$-analogue of $\{n\}$ which we will be considering is $\{n\}(q)=\{n\}_{s,t}(q)$  defined by $\{0\}(q)=0$, $\{1\}(q)=1$ and
\beq
\label{nq}
\{n\}(q)=s\{n-1\}(q) + tq^{n-2}\{n-2\}(q), 
\eeq
for $n\ge2$.  Cigler~\cite{cig:ncq,cig:qfp,cig:aam,cig:qpr}  introduced and studied these polynomials which have also been considered by Goyt and Sagan~\cite{gs:sps} and Goyt and Mathisen~\cite{gm:psq}.

To motivate the (generalized) Euler-Cassini identity, recall that a sequence of real numbers 
$( a_n)_{n\ge0}$ is called \emph{log concave} if it satisfies 
$$
a_{n}^{2} - a_{n-1}a_{n+1} \geq 0
$$ 
for all $n \geq 1$. Many sequences of combinatorial 
nature are log concave, for example any row of Pascal's triangle will do.  The identity 
$$
F_{rn}^{2} - F_{r(n+1)}F_{r(n-1)} = (-1)^{r} F_{r}^{2}
$$
shows that the sequence with  $a_{n} = F_{rn}$ is log concave for $r$ even.
See the articles of Brenti~\cite{bre:lus}, Stanley~\cite{sta:lus}, or Wilf~\cite{wil:gen}  for more details about log concavity and related issues. 

It is easy to see that, for sequences $(a_n)_{n\ge0}$ of positive reals, the log-concavity condition is equivalent to the seemingly stronger statement that
$$
a_n a_{n+m-1} \ge a_{n-1} a_{n+m}
$$
for all $m,n\ge1$.   The importance of the following generalization of the Euler-Cassini identity (which is the special case $s=t=q=1$) should now be clear.
\bth[\cite{cig:qfp}]
\label{EC}
We have

\eqqed{
\{n\}_{s,t}(q) \cdot \{n+m-1\}_{s,qt}(q)-\{n-1\}_{s,qt}(q) \cdot \{n+m\}_{s,t}(q)=
(-t)^{n-1}q^{\binom{n}2}\{m\}_{s,q^nt}(q).
}
\eth
\noindent  Note that as an immediate corollary, the sequence $\{n\}_{s,t}$ is log concave for all $t\le0$.

We are going to give a novel  proof for the $m=1$ case of this theorem using the Dodgson condensation technique for computing determinants.  So we will need a determinantal expression for $\{n\}(q)$.  We obtain this using a tri-diagonal matrix which is  a method common in the theory of special functions.
\bpr
\label{prop-det1}
The polynomial $\{ n \}(q)$ is given by 
$$
\{ n \}(q) = \det \begin{pmatrix}
s & -1 & 0 & \cdots & 0 & 0 \\
q t& s & -1 & \cdots & 0 & 0 \\
0  &  q^2 t & s & \cdots & 0 & 0 \\
\vdots  & \vdots & \vdots & \ddots & \vdots & \vdots \\
0 & 0 & 0 & \cdots & s & -1 \\
0 & 0 & 0 & \cdots &  q^{n-2}t& s 
\end{pmatrix} 
$$
where the tridiagonal matrix is of size $(n-1) \times (n-1)$.
\epr
\begin{proof}
By expansion about the last column, one easily verifies that the determinant satisfies the same initial conditions and recurrence as 
$\{ n \}(q)$.  
\end{proof}

For the final piece of background, we recall the method of Dodgson condensation. For any 
$n \times n$ matrix $A$, let $A_{r}(k,\ell)$ be the $r \times r$ connected
submatrix whose upper leftmost corner is the entry $a_{k,\ell}$.  If $\det A_{n-2}(2,2)\neq0$ then
\beq
\label{Dodgson}
\det A = \frac{\det A_{n-1}(1,1) \det A_{n-1}(2,2) - 
\det A_{n-1}(1,2) \det A_{n-1}(2,1)}{\det A_{n-2}(2,2)}. 
\eeq
Applications of this method can be found in the papers of Amdeberhan and Zeilberger~\cite{az:dtl} and Zeilberger~\cite{zei:rca}.

\noindent\emph{Proof of Theorem~\ref{EC} for $m=1$.}
To prove
\beq
\label{Euler-Cassini}
\{n\}_{s,t}(q) \cdot \{n\}_{s,qt}(q) -\{n-1\}_{s,qt}(q) \cdot \{n+1\}_{s,t}(q)=
(-t)^{n-1}q^{\binom{n}2},
\eeq
just apply equation~\ree{Dodgson} to the determinant for $\{n+1\}_{s,t}(q)$.  The result is
$$
\{n+1\}_{s,t}(q)=\frac{\{n\}_{s,t}(q) \cdot \{n\}_{s,qt}(q) - (-1)^{n-1}t^{n-1} q^{\binom{n}{2}}}{\{n-1\}_{s,qt}(q)}
$$
and applying a little algebra finishes the proof. 
\hqedm

It would be very interesting to prove the full version of Theorem~\ref{EC} in a similar manner.  This would perhaps require a more general version of condensation.

We now provide an application of the Euler-Cassini identity to infinite series.  We will need the following slight variant of Theorem~\ref{EC} when $q=1$.  It can be proved by adapting  Cigler's original proof.
\ble
\label{gEC}
Let $r,m,n\in\bbP$ and $s,t$ be arbitrary integers.  We have

\medskip

\eqqed{
\{rn\}\cdot \{r(n+m-1)\} -\{r(n-1)\} \cdot \{r(n+m)\} = (-t)^{r(n-1)} \{r\}\{rm\}.
}
\ele

Various authors have considered the following Fibonacci analogue of the Riemann zeta function
$$
\ze_F(z)=\sum_{k=0}^{\infty} \frac{1}{F_k^z}.
$$
See the article of Wu and Zhang~\cite{wz:sii} and references therein.  In particular, there has been interest in finding estimates of the tails of such series for positive integers $z$.  Holliday and Komatsu~\cite{hk:srg} considered what could be said for the Fibonacci polynomials where $t=1$ (recall~\ree{Fn(x)}) and proved the following result.
\bth[\cite{hk:srg}]
If $t=1$ and $s, n\in\bbP$, then
$$
\left\lfloor\left(\sum_{k=n}^{\infty}\frac1{\{k\}}_{s,1}\right)^{-1}\right\rfloor
=\{n\}_{s,1}-\{n-1\}_{s,1}-\de_O(n),
$$ 
and

\medskip

\eqqed{
\left\lfloor\left(\sum_{k=n}^{\infty}\frac1{\{k\}_{s,1}^2}\right)^{-1}\right\rfloor
=s\{n\}_{s,1}\{n-1\}_{s,1}-\de_E(n).}
\eth

Holliday and Komatsu also asked if their theorem could be generalized to other $t$ and we will do this for the first summation.  In addition, our results cover a more general class of sums and the proofs, based on Lemma~\ref{gEC}, will be much simpler than the ones given in~\cite{hk:srg}.  To see the equivalence of our second sum when $r=1$ with the one of Holliday and Komatsu, we note that if $t=1$ then
$$
\{n\}^2-\{n-1\}^2+(-1)^n=s\{n\}\{n-1\},
$$
an identity which is easily proved using Proposition~\ref{nXY} and equation~\ree{stXY}.

\bth
If $s\ge t\ge1$ and $n,r\in\bbP$ then
$$
\left\lfloor\left(\sum_{k=n}^{\infty}\frac1{\{rk\}}_{s,t}\right)^{-1}\right\rfloor
=\{rn\}_{s,t}-\{r(n-1)\}_{s,t}-\de_E(r(n-1)).
$$ 
If $t=1$ and $s,n,r\in\bbP$ then
$$
\left\lfloor\left(\sum_{k=n}^{\infty}\frac1{\{rk\}_{s,1}^2}\right)^{-1}\right\rfloor
=\{rn\}_{s,1}^2-\{r(n-1)\}_{s,1}^2-\de_E(r(n-1)).
$$ 
\eth
\bprf
First note that both series must converge by comparison with the known convergent series $\sum_{k\ge1}1/F_k$.  

We now consider the first series.  We will only give details for the case when $r(n-1)$ is even as the odd case is similar.  For ease of notation, let
$$
A(n)=\sum_{k=n}^{\infty}\frac1{\{rk\}}.
$$
It suffices to show that
$$
\{rn\}-\{r(n-1)\}-1\le\frac{1}{A(n)}<\{rn\}-\{r(n-1)\}.
$$
Note that since $s,t$ are positive, so are all the $\{rk\}$ and inequalities will not be affected when multiplying by them.

We first deal with the right-hand inequality.  Multiply through by $A(n)$ and then cancel the $1$ on the right with the first term of the series $\{rn\}A(n)$.  For $m\ge1$ we then compare the term for $k=n+m$ in $\{rn\}A(n)$  with the term for $k=n+m-1$ in $\{r(n-1)\}A(n)$ to see that it suffices to show 
\beq
\label{ineq1}
0<\frac{\{rn\}}{\{r(n+m)\}}-\frac{\{r(n-1)\}}{\{r(n+m-1)\}}.
\eeq
But this is true by Lemma~\ref{gEC} and the fact that $r(n-1)$ is even.

We now apply the same procedure as in the previous paragraph to the left-hand inequality and reduce it to proving
\beq
\label{ineq2}
\frac{\{rn\}}{\{r(n+m)\}}-\frac{\{r(n-1)\}}{\{r(n+m-1)\}}\le\frac{1}{\{r(n+m-1)\}}.
\eeq
Cross-multiplying and using Lemma~\ref{gEC} again as well as the parity of $r(n-1)$, we see that it suffices to prove 
$ t^{r(n-1)}\{r\}\{rm\}\le \{r(n+m)\}$.  Using Theorem~\ref{thm-recu1} and the fact that $t$ is positive gives
$$
\{r(n+m)\}\ge \{rn+1\}\{rm\} \ge \{r(n-1)+2\}\{r\}\{rm\} \ge t^{r(n-1)+1}\{r\}\{rm\}
$$
where the last inequality comes from the fact that, by Proposition~\ref{{n}-st} and $s\ge t\ge 1$, we have
$\{l\}\ge s^{l-1}\ge t^{l-1}$.  Thus we are done with the first series.

The proof for the second series has many similarities, so we will only mention the places where they differ.  Assume, again, that $r(n-1)$ is even.  To obtain the squared version of~\ree{ineq1}, merely move the negative fraction onto the other side of the inequality (which we know to be true from the first half of the proof) and square both sides.

When each fraction in inequality~\ree{ineq2} is replaced by its square one obtains, after clearing denominators,
$$
\{rn\}^2 \{r(n+m-1)\}^2 - \{r(n-1)\}^2 \{r(n+m)\}^2 \le \{rn+rm\}^2.
$$
Factoring the left-hand side and applying Lemma~\ref{gEC} once again, this time with $t=1$, gives the equivalent inequality
$$
\{r\} \{rm\} \cdot(\{rn\} \{r(n+m-1)\} + \{r(n-1)\} \{r(n+m)\})  \le \{rn+rm\}^2.
$$ 
It is easy to prove by induction that under the restrictions on $s,t$ the sequence $\{n\}$ is weakly increasing.   Using this observation as well as repeated application of Theorem~\ref{thm-recu1}, we obtain
\begin{align*}
\{r(n+m)\}^2
&=\{r(n+m)\}\cdot (\{r(n+m-1)\}\{r+1\}+\{r(n+m-1)-1\}\{r\})\\
&\ge \{rn\} \{rm\} \{r(n+m-1)\} \{r\} + \{r(n+m)\}\{r(n-1)\}\{rm\}\{r\}
\end{align*}
and factoring out $\{r\}\{rm\}$ leads to the desired conclusion.
\eprf

We believe that there are analogues of these results for other values of $s,t$.  The restriction $(s,t)\neq(2,-1)$ in the following conjecture is to ensure that the series converges.
\bcon
If $s> t\ge1$ with $(s,-t)\neq(2,-1)$ and $n,r\in\bbP$ then
$$
\left\lfloor\left(\sum_{k=n}^{\infty}\frac1{\{rk\}}_{s,-t}\right)^{-1}\right\rfloor
=\{rn\}_{s,-t}-\{r(n-1)\}_{s,-t}-1.
$$ 
If $t=-1$ and $s,n,r\in\bbP$ then
$$
\left\lfloor\left(\sum_{k=n}^{\infty}\frac{1}{\{rk\}_{s,-1}^2}\right)^{-1}\right\rfloor
=\{rn\}_{s,-1}^2-\{r(n-1)\}_{s,-1}^2-1.
$$ 
\econ

%
%

\section{Fundamental properties of $\binst{n}{k}$}
\label{sec-fpb}

We now return to the generalized Fibonomial coefficients.  We will describe various recursions which they satisfy as well as analogues of the binomial theorem and Chu-Vandermonde summation.

It is not clear from the definition that $\left\{ n \atop k \right\}$ is a polynomial in $s,t$ with nonnegative integral coefficients.  This will follow by an easy induction using the first recursion in Theorem~\ref{rrs}.  The Lucanomials are not so well behaved.   For instance, 
$$
\left\langle  4 \atop 2 \right\rangle = 
\frac{(s^{2}+3t)(s^{4}+4s^{2}t+2t^{2})}{s^{2}+2t}.
$$
can not be brought to  polynomial form.  Also, we could introduce a $q$-analogue $\binst{n}{k}(q)$ of the generalized Fibonomials by using the the polynomials $\{n\}(q)$ defined by~\ree{nq}.  But then it is easy to check that $\binst{6}{3}(q)$ is not a polynomial.  This is one of the reasons we have decided to mainly consider the case $q=1$ in this work.

Our first property of generalized Fibonomials is their symmetry.  This follows immediately from their definition.
\bpr
If $0 \leq k \leq n$ then

\medskip

\eqqed{
\left\{ n \atop {k}  \right\} = \left\{ n \atop {n-k} \right\}.
}
\epr

Next we consider two recursions for the Fibonomials.
\bth
\label{rrs}
For $m,n\ge1$ we have
\begin{align}
\label{rr1}
\left\{ {m+n} \atop m \right\}
& =   \{ n+1 \} \left\{ {n+m-1} \atop {m-1} \right\}  + t \{ m-1 \}  \left\{ {m+n-1} \atop {m} \right\} 
\\[10pt]
\label{rr2}
&= Y^n \left\{ {n+m-1} \atop {m-1} \right\} + X^m  \left\{ {m+n-1} \atop {m} \right\}.
\end{align}
In particular they are polynomials in $s$ and $t$. 
\eth
\begin{proof}
The first recursion follows easily from Theorem \ref{thm-recu1}.  The second can be obtained from the $q$-binomial recursion
$$
\gaus{m+n}{m}=\gaus{m+n-1}{m-1}+q^m\gaus{m+n-1}{m}
$$
via the substitution~\ree{gau}.
\end{proof}

The next result gives two analogues of the binomial theorem. 
\bth
Letting $z$ be an indeterminate, the genearlized Fibonomials satisfy
\beq
\label{bt}
(1+X^{n-1} z) (1+ X^{n-2} Y z)\cdots (1+Y^{n-1}z)=\sum_{k=0}^n  (-t)^{\binom{k}{2}} \binst{n}{k}  z^k
\eeq
and
\beq
\label{negbt}
\frac{1}{(1-X^{n-1} z) (1- X^{n-2} Y z)\cdots (1-Y^{n-1}z)}=\sum_{k=0}^\infty \binst{n+k-1}{k}  z^k.
\eeq
\eth
\bprf
We will indicate how to prove the first recursion as the second is similar.  One approach is to use~\ree{rr2} and induction.  Alternatively, one can start with the $q$-binomial theorem in the form
$$
(1+z)(1+qz)\cdots(1+q^{n-1}z)=\sum_{k=0}^n q^{\binom{k}{2}} \gaus{n}{k}  z^k,
$$
substitute $Y^{n-1}z$ for $z$, use~\ree{gau}, and clear denominators.  The details of both proofs are routine and so these are left to the reader.
\eprf

It is interesting to note that when $t=-1$ we can write~\ree{negbt} using generalized Fibonomials with negative upper indices.  
We note that in this case, the generalized Fibonacci polynomial sequence is  called an $\ell$-sequence (where $s=\ell$) and such sequences are intimately related with lecture hall partitions; see the papers by Bousquet-M\'elou and  Eriksson~\cite{bme:lhp1,bme:lhp2} and by Savage and Yee~\cite{sy:ept} for more information.
We first extend the sequence $\{n\}$ to negative integers by insisting that the recursion~\ree{nst-rr} continues to hold.  In this case, it is not hard to show by induction that, for $n\ge0$, 
$$
\{-n\}=\frac{-\{n\}}{(-t)^n}.
$$
So when $t=-1$ the $\{-n\}$ are polynomials in $s,t$ with integral coefficients.  In this case, let
$$
\binst{-n}{k}_{s,-1} \stackrel{\rm}{=} \frac{\{-n\}\{-n-1\}\cdots\{-n-k+1\}}{\{k\}!}=(-1)^k\binst{n+k-1}{k}_{s,-1}.
$$
Hence~\ree{negbt} becomes
$$
\frac{1}{(1-X^{n-1} z) (1- X^{n-2} Y z)\cdots (1-Y^{n-1}z)}=\sum_{k=0}^\infty \binst{-n}{k}_{s,-1}  (-z)^k.
$$

The specialization $t=-1$ also permits us to obtain nice analogues of the formulas for the sum and alternating sum of a row of Pascal's triangle.  This is because, by~\ree{stXY}, we have $XY=1$ which permits simplifications.  Note also that if $s=2$ and $t=-1$ then it is easy to prove that $\langle n\rangle_{2,-1}=2$ for all $n\ge0$, so that in this case the following identities reduce to the usual ones.
\bco
If $t=-1$ then
$$
\sum_{k=0}^n \binst{n}{k}_{s,-1}=(1+\de_O(n))\prod_{i=1}^{\lf n/2\rf}(2+\langle n-2i+1 \rangle)
$$
and
$$
\sum_{k=0}^n (-1)^k \binst{n}{k}_{s,-1}=\de_E(n)\prod_{i=1}^{\lf n/2\rf}(2-\langle n-2i+1 \rangle).
$$
\eco
\bprf
As usual, we just prove the first identity.    Letting $t=-1$ and $z=1$ in equation~\ree{bt} gives
$$
\sum_{k=0}^n  \binst{n}{k}_{s,-1}  = (1+X^{n-1} ) (1+ X^{n-2} Y )\cdots (1+Y^{n-1})
$$
If $i<j$ then we have, using Proposition~\ref{nXY},
$$
(1+X^i Y^j)(1+X^j Y^i) = 1 + X^i Y^i(X^{j-i}+Y^{j-i}) + X^{i+j} Y^{i+j} = 2 + \langle j-i \rangle.
$$
Pairing up such factors and remembering that there will be an unpaired factor when $n$ is odd completes the proof.
\eprf

Returning to an arbitrary $t$, we can use any identity for $q$-integers and $q$-binomials to derive a corresponding one for $\{n\}$ and $\binst{n}{k}$.  For example, the $q$-Chu-Vandermonde summation
$$
\gaus{m+n}{k} = \sum_i  q^{i(m-k+i)}\gaus{m}{k-i} \gaus{n}{i}
$$
gives rise to the following result.
\bth
We have

\medskip

\eqqed{
\binst{m+n}{k}=\sum_i (-t)^{i(i-k)} X^{mi} Y^{n(k-i)} \binst{m}{k-i} \binst{n}{i}.
}
\eth

We end this section by mentioning that Proposition~\ref{gcd} can be used to prove various divisibility properties of the Fibonomials $\binst{n}{k}$.  An example is the following primality testing condition.  A proof for ordinary binomial coefficients can be found in the article of Dilcher and Stolarsky~\cite{ds:ptt}.  Because of Proposition~\ref{gcd}, the demonstration  carries over mutatis mutandis to the general case. 
\bpr
Consider $s,t$ as variables.  The positive integer $p$ is prime if and only if $\{p\}$ divides $\binst{p}{k}$ for each $0<k<p$.\hqed
\epr

%
%

\section{Catalan numbers}
\label{sec-Catalan}

In this section, we will consider an $s,t$-analogue of the Catalan numbers suggested by Louis Shapiro~[private communication].  
We will then investigate the $2$-adic valuation of these generalized Catalan numbers for various values of $s$ and $t$, extending a well-known theorem in the case $s=2, t=-1$.

Recall that the  \emph{Catalan numbers} are defined by
$$
C_{n} = \frac{1}{n+1} \binom{2n}{n}
$$
for $n\ge0$.  By analogy, define the \emph{generalized Catalan numbers} to be
$$
C_{\{n\}} = \frac{1}{\{ n+1 \}} \left\{ {2n} \atop {n} \right\}.
$$
 Ekhad~\cite{ekh:ssl} was the first to note that these must be polynomials in $s,t$ with nonnegative integral coefficients, because setting $m=n$ in Theorem~\ref{thm-recu1} and doing some algebraic manipulation shows that
$$
C_{\{ n \}} = \left\{ {2n-1} \atop {n-1} \right\} + t 
\left\{ {2n-1} \atop {n-2} \right\}.
$$
So we can ask about the arithmetic properties of $C_{\{n\}}$ for integers $s,t$.

To state and prove our results, we will need some notation.  First of all, let
$$
\ze_\bb(n)=\text{the number of nonzero digits of $n$ in the base $\bb$}.
$$
Also, let
$$
\ka_\bb(m+n)=\text{the number of carries in doing the addition $(m)_\bb+(n)_\bb$}.
$$

It will also be useful to introduce the expansion of $n$ in an unusual base for the integers.  Given any infinite increasing sequence of positive integers $\bb=(b_0,b_1,\ldots)$ such that $b_0=1$ and $b_i|b_{i+1}$ for $i\ge1$ we consider the \emph{expansion of $n$ in base $\bb$} to be $n=\sum_{i\ge0} n_i b_i$ where $0\le n_i< b_{i+1}/b_i$ for all $i\ge0$.  We will utilize the shorthand $(n)_\bb =(n_0,n_1,\ldots)$ for the digits $n_i$ in this expansion.  If $\bb$ consists of powers of an integer $m>1$ then we will merely write $(n)_m$.  And we will omit the subscript entirely if the base is clear from context.   A particular base of interest to us is the one where 
$$
\bF=(1,3,3\cdot 2, 3\cdot 2^2, 3\cdot 2^3,\ldots).
$$
This base arises naturally from studying the fractal nature of the Fibonomial triangle modulo $2$.  See the paper of Chen and Sagan~\cite{cs:fnf} for details.   Other number-theoretic functions which depend on choosing a base will follow the same conventions.  

The next result is well-known.  See the paper by Deutsch and Sagan~\cite{ds:ccm} for a (mostly) combinatorial proof.
\bth
\label{thm-nu_2(C_n)}
If $C_n$ is an ordinary Catalan number then

\eqqed{
\nu_2(C_n)=\ze_2(n+1)-1.
}
\eth

We now prove analogues of this theorem for $C_{\{n\}}$ using results from Section~\ref{sec-arithmetic}.  We  first recall a  famous theorem of Kummer~\cite{kum:ear} which is needed in the sequel.
\bth
\label{thm-Kummer}
For any prime, $p$, we have

\eqqed{
\nu_p\left(\ \binom{m+n}{n}\ \right) = \ka_p(m+n).
}
\eth

We start  by considering what happens if $s$ and $t$ are both odd.

\bth
Let $s$ and $t$ be odd. We have
$$
\nu_2(C_{\{n\}})=
\begin{cases}
\ze_\bF(n+1)+\nu_2(\{6\})-3 & \text{if $t\equiv 1\ (\Mod 4)$ and $n\equiv 3, 4\ (\Mod 6)$,}\\
\ze_\bF(n+1)-1 & \text{else.}
\end{cases}
$$
\eth
\bprf
Since
$$
\nu_2(C_{\{n\}})=\nu_2(\{2n\}!)-\nu_2(\{n\}!)-\nu_2(\{n+1\}!),
$$
we can apply Corollary~\ref{cor-sodd-todd}.  As usual, we will just supply details when $t\equiv 3\ (\Mod 4)$.  In this case, the terms from the corollary containing a factor of $\nu_2(\{3\})$ cancel in the above equation. As a result, 
$$
\nu_2(C_{\{n\}})=\nu_2\left(\ \binom{\lf 2n/3 \rf}{\lf n/3\rf}\ \right)=
\begin{cases}
\ze_2( \lf n/3\rf )&\text{if $n\equiv 0, 1\ (\Mod 3)$,}\\
\ze_2( \lf n/3\rf   + 1)-1&\text{if $n\equiv 2\ (\Mod 3)$,}
\end{cases}
$$
where the second equality comes from Kummer's Theorem and the fact that when adding $k$ to itself in base two the number of carries is the number of nonzero digits.  We now translate from base $2$ to base $\bF$.  We must consider the congruence classes modulo three individually. We will do $n\equiv 0\ (\Mod 3)$ and leave the rest to the reader.  In this case
$n=3k$ where $(k)_2=(k_0,k_1,\ldots)$.  So $(n+1)_\bF =(3k+1)_\bF=(1,k_0,k_1,\ldots)$.  Thus 
$$
\ze_2( \lf n/3\rf )=\ze_2(k)=\ze_\bF(3k+1)-1=\ze_\bF(n+1)-1
$$
which is what we wished to prove.
\eprf

Now we consider $s$ and $t$ of opposite parity.   As remarked before, if $s$ is odd and $t$ is even then  $C_{\{n\}}$ is 
always odd.  If $s$  is even and $t$ is odd then Corollary~\ref{cor-seven-todd} and Theorem~\ref{thm-nu_2(C_n)} also make evaluation a simple matter.  Thus we arrive at the following theorem.

\bth
Let $s$ and $t$ be of opposite parity.  We have

\smallskip

\eqqed{
\nu_2(C_{\{n\}})=\case{\ze_2(n+1)-1}{if $t$ is odd,}{0}{if $t$ is even.}
}
\eth

We conclude with a very interesting question that Shapiro raised when he defined the Catalans $C_{\{n\}}$.  As mentioned earlier, there is a nice combinatorial interpretation of $\binst{n}{k}$ using tilings~\cite{ss:cib}.  But it remains an open problem to find a combinatorial interpretation for $C_{\{n\}}$.  This is especially puzzling given the plethora of combinatorial interpretations for the ordinary Catalan numbers.


\begin{thebibliography}{10}

\bibitem{az:dtl}
Tewodros Amdeberhan and Doron Zeilberger.
\newblock Determinants through the looking glass.
\newblock {\em Adv. in Appl. Math.}, 27(2-3):225--230, 2001.
\newblock Special issue in honor of Dominique Foata's 65th birthday
  (Philadelphia, PA, 2000).

\bibitem{bme:lhp1}
Mireille Bousquet-M\'elou and Kimmo Eriksson.
\newblock Lecture hall partitions.
\newblock {\em Ramanujan J.}, 1(1):101--111, 1997.

\bibitem{bme:lhp2}
Mireille Bousquet-M\'elou and Kimmo Eriksson.
\newblock Lecture hall partitions {II.}
\newblock {\em Ramanujan J.}, 1(2):165--185, 1997.

\bibitem{bre:lus}
Francesco Brenti.
\newblock Log-concave and unimodal sequences in algebra, combinatorics, and
  geometry: an update.
\newblock In {\em Jerusalem combinatorics '93}, volume 178 of {\em Contemp.
  Math.}, pages 71--89. Amer. Math. Soc., Providence, RI, 1994.

\bibitem{cs:fnf}
Xi~Chen and Bruce~E. Sagan.
\newblock On the fractal nature of {F}ibonomial coefficients.
\newblock In preparation.

\bibitem{cig:ncq}
Johann Cigler.
\newblock A new class of {$q$}-{F}ibonacci polynomials.
\newblock {\em Electron. J. Combin.}, 10:Research Paper 19, 15 pp., 2003.

\bibitem{cig:qfp}
Johann Cigler.
\newblock {$q$}-{F}ibonacci polynomials.
\newblock {\em Fibonacci Quart.}, 41(1):31--40, 2003.

\bibitem{cig:aam}
Johann Cigler.
\newblock Some algebraic aspects of {M}orse code sequences.
\newblock {\em Discrete Math. Theor. Comput. Sci.}, 6(1):55--68 (electronic),
  2003.

\bibitem{cig:qpr}
Johann Cigler.
\newblock {$q$}-{F}ibonacci polynomials and the {R}ogers-{R}amanujan
  identities.
\newblock {\em Ann. Comb.}, 8(3):269--285, 2004.

\bibitem{ds:ccm}
Emeric Deutsch and Bruce~E. Sagan.
\newblock Congruences for {C}atalan and {M}otzkin numbers and related
  sequences.
\newblock {\em J. Number Theory}, 117(1):191--215, 2006.

\bibitem{ds:ptt}
Karl Dilcher and Kenneth~B. Stolarsky.
\newblock A {P}ascal-type triangle characterizing twin primes.
\newblock {\em Amer. Math. Monthly}, 112(8):673--681, 2005.

\bibitem{dod:cd}
C.~L. Dodgson.
\newblock Condensation of determinants, being a new and brief method for
  computing their arithmetic values.
\newblock {\em Proc. Royal Soc. of London}, 15:150--155, 1866.

\bibitem{ekh:ssl}
Shalosh Ekhad.
\newblock
  The~{S}agan-{S}avage~{L}ucas-{C}atalan~polynomials~have~positive~coefficient%
s.
\newblock Preprint
  {\texttt{http://www.math.rutgers.edu/\raisebox{-2pt}{$\tilde{\rule{5pt}{0pt}}$}zeilberg/mamarim/mamarimhtml/bruce.htm%
l}}.

\bibitem{gm:psq}
Adam~M. Goyt and David Mathisen.
\newblock Permutation statistics and {$q$}-{F}ibonacci numbers.
\newblock {\em Electron. J. Combin.}, 16(1):Research Paper 101, 15 pp., 2009.

\bibitem{gs:sps}
Adam~M. Goyt and Bruce~E. Sagan.
\newblock Set partition statistics and {$q$}-{F}ibonacci numbers.
\newblock {\em European J. Combin.}, 30(1):230--245, 2009.

\bibitem{hl:dpg}
Verner~E. Hoggatt, Jr. and Calvin~T. Long.
\newblock Divisibility properties of generalized {F}ibonacci polynomials.
\newblock {\em Fibonacci Quart.}, 12:113--120, 1974.

\bibitem{hk:srg}
Sarah~H. Holliday and Takao Komatsu.
\newblock On the sum of reciprocal generalized {F}ibonacci numbers.
\newblock {\em Integers}, 11(4):441--455, 2011.

\bibitem{kos:fln}
Thomas Koshy.
\newblock {\em Fibonacci and {L}ucas numbers with applications}.
\newblock Pure and Applied Mathematics (New York). Wiley-Interscience, New
  York, 2001.

\bibitem{kum:ear}
E.~E. Kummer.
\newblock {\"U}ber die {E}rg\"anzungss\"atze zu den allgemeinen
  reciprocit\"atsgesetzen.
\newblock {\em J. reine angew. Math.}, 44:93--146, 1852.

\bibitem{luc:tfn1}
Edouard Lucas.
\newblock Theorie des {F}onctions {N}umeriques {S}implement {P}eriodiques.
\newblock {\em Amer. J. Math.}, 1(2):184--196, 1878.

\bibitem{luc:tfn2}
Edouard Lucas.
\newblock Theorie des {F}onctions {N}umeriques {S}implement {P}eriodiques.
  [{C}ontinued].
\newblock {\em Amer. J. Math.}, 1(3):197--240, 1878.

\bibitem{luc:tfn3}
Edouard Lucas.
\newblock Theorie des {F}onctions {N}umeriques {S}implement {P}eriodiques.
\newblock {\em Amer. J. Math.}, 1(4):289--321, 1878.


\bibitem{mol:nfs}
Victor~H. Moll.
\newblock {\em Numbers and functions}.
\newblock Student Mathematical Library. American Mathematical Society,
  Providence, RI, 2012.
\newblock Special Functions for Undergraduates.

\bibitem{ss:cib}
Bruce~E. Sagan and Carla~D. Savage.
\newblock Combinatorial interpretations of binomial coefficient analogues
  related to {L}ucas sequences.
\newblock {\em Integers}, 10:A52, 697--703, 2010.

\bibitem{sy:ept}
Carla~D. Savage and Ae~Ja Yee.
\newblock Euler's partition theorem and the combinatorics of
  {$\ell$}-sequences.
\newblock {\em J. Combin. Theory Ser. A}, 115(6):967--996, 2008.

\bibitem{sta:lus}
Richard~P. Stanley.
\newblock Log-concave and unimodal sequences in algebra, combinatorics, and
  geometry.
\newblock In {\em Graph theory and its applications: East and West (Jinan,
  1986)}, volume 576 of {\em Ann. New York Acad. Sci.}, pages 500--535. New
  York Acad. Sci., New York, 1989.

\bibitem{wil:gen}
Herbert~S. Wilf.
\newblock {\em generatingfunctionology}.
\newblock Academic Press Inc., Boston, MA, second edition, 1994.

\bibitem{wz:sii}
Zhengang Wu and Wenpeng Zhang.
\newblock Several identities involving the {F}ibonacci polynomials and {L}ucas
  polynomials.
\newblock {\em J. Inequal. Appl.}, 205, 14 pp., 2013.


\bibitem{zei:rca}
Doron Zeilberger.
\newblock Reverend {C}harles to the aid of {M}ajor {P}ercy and {F}ields
  medalist {E}nrico.
\newblock {\em Amer. Math. Monthly}, 103(6):501--502, 1996.

\end{thebibliography}
\end{document}